\documentclass[11pt,fleqn]{article}

\usepackage[cp1251]{inputenc}

\usepackage{amsmath,amssymb,amsthm,enumerate,epsfig,graphics,cite}
\usepackage[ukrainian,english]{babel}
\usepackage{amsfonts}

\newcommand{\todo}[1][\null]{\ensuremath{\clubsuit}}

\newcommand{\noprint}[1]{}

\begin{document}
\begin{center}
\Large\bf
Two-dimensional generalized gamma function and its applications
\end{center}
\begin{center}
\it
Artem M. Ponomarenko
\end{center}

\begin{center}
National Technical University of Ukraine
"Igor Sikorsky Kyiv Polytechnic Institute"
\end{center}

{\it

In this article, we present a new two-dimensional generalization of the gamma function based on the product of the one-dimensional generalized beta function and the one-dimensional generalized gamma function.
As will become clear later, this generalization is also a generalization
the famous formula that gives the connection between the classical gamma and beta functions. Next we present
properties of this generalization, some series for the generalized beta function. As a practical application of the two-dimensional generalized gamma function, we will show how it
can be used to represent a fairly wide class of double integrals in the form of functional series.

}
\medskip


\medskip

The properties and applications of classical special functions are well studied and are presented, for example, in \cite {1}, \cite {2}, \cite {3}, \cite {4}. Among other things, this literature presents many properties of the classical special functions $\Gamma(\alpha)$ and $B(\alpha,\beta)$ and their various generalizations. In this article we will consider a new generalization of the Gamma function. We will also obtain a generalization of the formula for the relationship between the classical gamma function and the beta function \cite {3} $\Gamma(\alpha)\Gamma(\beta)=B(\alpha,\beta)\Gamma(\alpha+\beta)$ to the case of any integrand functions.

\medskip

Throughout the article we will consider the parameters $\alpha$, $\beta$, $\gamma$, $\omega$, $\nu$, $\lambda$, $a$, $b$ as real numbers only, with the exception of Remark 5 and the final application for hypergeometric functions at the end of the article.

\medskip
\section{Introduction of a two-dimensional generalization of the gamma function and study of its properties.}
We will consider a one-dimensional generalization of the gamma function in the form
\begin{gather}\label{eq1dd3d3}
\Gamma_{g(\cdot)}(\omega):=\int_{0}^{\infty}g(x)x^{\omega-1}e^{-x}dx, \quad \omega>0,
\end{gather}
one-dimensional generalization of the beta function in the form
\begin{gather}\label{eq2}
B_{f(\cdot)}(\alpha,\beta):=\int_{0}^{1}f(x)x^{\alpha-1}(1-x)^{\beta-1}dx, \quad \alpha>0, \quad \beta>0,
\end{gather}
and two-dimensional generalization of the Gamma function
\begin{gather}\label{eq1dd3d4}
\Gamma_{\mathbb{R}^{2}_{+}; \Omega \left(\ast_{_{\nu}}, \cdot_{_{\omega}}\right)}(\nu, \omega; \lambda):=\overset {} {\underset {\{T=(x,y)\in \mathbb{R}^{2} \hspace{1mm}  : \hspace{1mm} x\geq 0, y\geq 0 \} }{\int\int} }  \Omega (y,x) y^{\nu-1}x^{\omega-1}(x+y)^{\lambda}e^{-x-y}dydx, \quad \nu>0,
\end{gather}
$$
\omega>0, \quad \lambda\geq 0.
$$
Obviously, if $\Omega\equiv 1$ and $\gamma=0$, then our generalized Gamma-function $\Gamma_{\mathbb{R}^{2}_{+};}(\nu, \omega; 1)$ is transformed into the product of two classical Gamma-functions: $\Gamma(\alpha)\Gamma(\beta)$.

In this article we will set the task of finding the product of generalizations \eqref{eq1dd3d3} and \eqref{eq2}
 $$
 B_{f(\cdot)}(\alpha,\beta)\Gamma_{g(\cdot)}(\omega)
 $$
 in the form \eqref{eq1dd3d4}. That is, we will find for the function $\Omega$ such a condition that this product looks like a special case of the right-hand side of the formula \eqref{eq1dd3d4}. Thus, we obtain a generalization of the well-known formula for the relationship between the classical gamma function and the beta function $\Gamma(\alpha)\Gamma(\beta)=B(\alpha,\beta)\Gamma(\alpha+\beta)$.

To do this, we will consider the following generalizations of the gamma function of the form
\begin{gather}\label{eq1}
\Gamma_{g(\cdot)}(\alpha+\beta+\gamma):=\int_{0}^{\infty}g(x)x^{\alpha+\beta+\gamma-1}e^{-x}dx, \quad \alpha>0, \quad \beta>0, \quad \gamma\geq 0,
\end{gather}
and provided that integrals \eqref{eq1}, \eqref{eq2} exist, and we will obtain a double integral, which will be equal to their product. We will also consider integrals
$$
I_{\mathbb{R}^{2}_{+}; f\left(\frac{\ast_{_{\alpha}}}{\cdot_{_{\beta}}+\ast_{_{\alpha}}}\right)g\left( \cdot_{_{\beta}}+\ast_{_{\alpha}} \right)}(\alpha,\beta;\gamma):=
$$
\begin{gather}\label{eq3d1}
:= \overset {} {\underset {\{T=(x,y)\in \mathbb{R}^{2} \hspace{1mm}  : \hspace{1mm} x\geq 0, y\geq 0 \} }{\int\int} }  f\left(\frac{y}{x+y}\right)g\left( x+y \right)y^{\alpha-1}x^{\beta-1}\left( x+y \right)^{\gamma}e^{-x-y}dydx,
\end{gather}
$$
I^{[1]}_{\mathbb{R}^{2}_{+}; f\left(\frac{\ast_{_{\alpha}}}{\cdot_{_{\beta}}+\ast_{_{\alpha}}}\right)g\left( \cdot_{_{\beta}}+\ast_{_{\alpha}} \right)}(\alpha,\beta;\gamma):=
$$
\begin{gather}\label{eq3}
:=\int_{0}^{\infty}\int_{0}^{\infty}f\left(\frac{y}{x+y}\right)g\left( x+y \right)y^{\alpha-1}x^{\beta-1}\left( x+y \right)^{\gamma}e^{-x-y}dydx,
\end{gather}
and
$$
I^{[2]}_{\mathbb{R}^{2}_{+}; f\left(\frac{\ast_{_{\alpha}}}{\cdot_{_{\beta}}+\ast_{_{\alpha}}}\right)g\left( \cdot_{_{\beta}}+\ast_{_{\alpha}} \right)}(\alpha,\beta;\gamma):=
$$
\begin{gather}\label{eq4}
:=\int_{0}^{\infty}\int_{0}^{\infty}f\left(\frac{y}{x+y}\right)g\left( x+y \right)y^{\alpha-1}x^{\beta-1}\left( x+y \right)^{\gamma}e^{-x-y}dxdy,
\end{gather}
where $\alpha>0$, $\beta>0$, $\gamma\geq 0$ for all three previous integrals.

We will see later in our study that in the generalization \eqref{eq1dd3d4} it is convenient to choose the order of variables $\Omega (y,x)$, and not $\Omega (x,y)$.

Let us introduce the following notation
\begin{gather}\label{eq4ddd1d1}
Q(x,y):=f\left(\frac{y}{x+y}\right)g\left( x+y \right).
\end{gather}

{\it Theorem 1.} {\it Let us assume that $Q(x,y)$ is a given function satisfying the following conditions:

1) the function} $Q(x,y)$ {\it exists, is non-negative, continuous at all points of the region } $\{T=(x,y)\in \mathbb{R}^{2} \hspace{1mm}  : \hspace{1mm} x\geq 0, y\geq 0 \} $, {\it with the possible exception of a finite number of points } $M_{1}$, $M_{2}$, ... , $M_{\kappa}$, $\kappa \in \mathbb{N}$, {\it in this region. At points} $M_{1}$, $M_{2}$, ... , $M_{\kappa}$ {\it the function} $Q(x,y)$ {\it tends to} $+\infty$;

{\it 2) integrals} \eqref{eq1}, \eqref{eq2}, \eqref{eq3d1}, \eqref{eq3}, \eqref{eq4} {\it exist and converge for all values} $\alpha>0$, $\beta>0$, $\gamma\geq 0$;

{\it 3) the following equality holds}
$$
I_{\mathbb{R}^{2}_{+}; f\left(\frac{\ast_{_{\alpha}}}{\cdot_{_{\beta}}+\ast_{_{\alpha}}}\right)g\left( \cdot_{_{\beta}}+\ast_{_{\alpha}} \right)}(\alpha,\beta;\gamma)=
$$
\begin{gather}\label{eq5}
=I^{[1]}_{\mathbb{R}^{2}_{+}; f\left(\frac{\ast_{_{\alpha}}}{\cdot_{_{\beta}}+\ast_{_{\alpha}}}\right)g\left( \cdot_{_{\beta}}+\ast_{_{\alpha}} \right)}(\alpha,\beta;\gamma)=I^{[2]}_{\mathbb{R}^{2}_{+}; f\left(\frac{\ast_{_{\alpha}}}{\cdot_{_{\beta}}+\ast_{_{\alpha}}}\right)g\left( \cdot_{_{\beta}}+\ast_{_{\alpha}} \right)}(\alpha,\beta;\gamma),
\end{gather}
{\it where} $\alpha>0$, $\beta>0$, $\gamma\geq 0$.

{\it Then the following formula will be valid }
\begin{gather}\label{eq6}
B_{f(\cdot)}(\alpha,\beta)\Gamma_{g(\cdot)}(\alpha+\beta+\gamma)=I^{[1]}_{\mathbb{R}^{2}_{+}; f\left(\frac{\ast_{_{\alpha}}}{\cdot_{_{\beta}}+\ast_{_{\alpha}}}\right)g\left( \cdot_{_{\beta}}+\ast_{_{\alpha}} \right)}(\alpha,\beta;\gamma),\quad \alpha>0,\hspace{1mm} \beta>0,\hspace{1mm} \gamma\geq 0.
\end{gather}

{\it Proof.}
$$
B_{f(\cdot)}(\alpha,\beta)=\int_{0}^{1}f\left( x \right)x^{\alpha-1}(1-x)^{\beta-1}dx=\int_{0}^{1}f\left( \frac{\Gamma(x+1)}{\Gamma(x)} \right)x^{\alpha-1}(1-x)^{\beta-1}dx=
$$
$$
=\int_{0}^{1}f\left( \frac{\int_{0}^{\infty}t^{x+1-1}e^{-t}dt}{\int_{0}^{\infty}\tau^{x-1}e^{-\tau}d\tau} \right)x^{\alpha-1}(1-x)^{\beta-1}dx=
$$
$$
=2\int_{0}^{1}f\left( \frac{\int_{0}^{\infty}t^{x^{2}+1-1}e^{-t}dt}{\int_{0}^{\infty}\tau^{x^{2}-1}e^{-\tau}d\tau} \right)x^{2\alpha-1}(1-x^{2})^{\beta-1}dx=
$$
\begin{gather}\label{eq7}
=2\int_{0}^{\frac{\pi}{2}}f\left( \frac{\int_{0}^{\infty}t^{\sin^{2}\varphi}e^{-t}dt}{\int_{0}^{\infty}\tau^{-\cos^{2}\varphi}e^{-\tau}d\tau} \right)\sin^{2\alpha-1}(\varphi) \cos^{2\beta-1}(\varphi) d\varphi.
\end{gather}
Let's multiply both sides of formula \eqref{eq7} by $g\left( r^{2} \right)r^{2\alpha+2\beta+2\gamma-1}e^{-r^{2}}$ and integrate it from $0$ to $\infty$. Then we obtain
$$
B_{f(\cdot)}(\alpha,\beta)\int_{0}^{\infty}g\left( r^{2} \right)r^{2\alpha+2\beta+2\gamma-1}e^{-r^{2}}dr=
$$
$$
=2\int_{0}^{\frac{\pi}{2}}\int_{0}^{\infty} f\left( \frac{\int_{0}^{\infty}t^{\sin^{2}\varphi}e^{-t}dt}{\int_{0}^{\infty}\tau^{-\cos^{2}\varphi}e^{-\tau}d\tau} \right) g\left( r^{2} \right) (r \sin\varphi)^{2\alpha-1} (r \cos\varphi)^{2\beta-1} r^{2\gamma} e^{-r^{2}} r dr d\varphi=
$$
$$
=2\int_{0}^{\frac{\pi}{2}}\int_{0}^{\infty} f\left( \frac{r^{2}\int_{0}^{\infty}t^{r^{2}\sin^{2}\varphi}e^{-t^{r^{2}}}t^{r^{2}-1}dt}{r^{2}\int_{0}^{\infty}\tau^{-r^{2}\cos^{2}\varphi}e^{-\tau^{r^{2}}}\tau^{r^{2}-1}d\tau} \right) g\left( r^{2} \right) (r \sin\varphi)^{2\alpha-1} (r \cos\varphi)^{2\beta-1} r^{2\gamma} e^{-r^{2}} r dr d\varphi=
$$
$$
=2\int_{0}^{\frac{\pi}{2}}\int_{0}^{\infty} f\left( \frac{((r \sin\varphi)^{2}+(r \cos\varphi)^{2})\int_{0}^{\infty}t^{r^{2}\sin^{2}\varphi}e^{-t^{(r \sin\varphi)^{2}+(r \cos\varphi)^{2}}}t^{(r \sin\varphi)^{2}+(r \cos\varphi)^{2}-1}dt}{((r \sin\varphi)^{2}+(r \cos\varphi)^{2})\int_{0}^{\infty}\tau^{-r^{2}\cos^{2}\varphi}e^{-\tau^{(r \sin\varphi)^{2}+(r \cos\varphi)^{2}}}\tau^{(r \sin\varphi)^{2}+(r \cos\varphi)^{2}-1}d\tau} \right)  \cdot
$$
$$
\cdot g\left( (r \sin\varphi)^{2}+(r \cos\varphi)^{2} \right)(r \sin\varphi)^{2\alpha-1} (r \cos\varphi)^{2\beta-1}\left((r \sin\varphi)^{2}+(r \cos\varphi)^{2}\right)^{\gamma}\cdot
$$
$$
\cdot e^{-\left((r \sin\varphi)^{2}+(r \cos\varphi)^{2}\right)} r dr d\varphi=
$$
$$
=2\int_{0}^{\infty}\int_{0}^{\infty} f\left( \frac{(x^{2}+y^{2})\int_{0}^{\infty}t^{y^{2}}e^{-t^{x^{2}+y^{2}}}t^{x^{2}+y^{2}-1}dt}{(x^{2}+y^{2})\int_{0}^{\infty}\tau^{-x^{2}}e^{-\tau^{x^{2}+y^{2}}}\tau^{x^{2}+y^{2}-1}d\tau} \right) g\left( x^{2}+y^{2} \right) \cdot
$$
$$
\cdot y^{2\alpha-1} x^{2\beta-1} \left(x^{2}+y^{2}\right)^{\gamma} e^{-x^{2}-y^{2}}  dy dx =
$$
$$
=\frac{1}{2}\int_{0}^{\infty}\int_{0}^{\infty} f\left( \frac{(x+y)\int_{0}^{\infty}t^{y}e^{-t^{x+y}}t^{x+y-1}dt}{(x+y)\int_{0}^{\infty}\tau^{-x}e^{-\tau^{x+y}}\tau^{x+y-1}d\tau} \right) g\left( x+y \right) y^{\alpha-1} x^{\beta-1} (x+y)^{\gamma} e^{-x-y}  dy dx=
$$
$$
=\frac{1}{2}\int_{0}^{\infty}\int_{0}^{\infty} f\left( \frac{\int_{0}^{\infty}t^{\frac{y}{x+y}}e^{-t}dt}{\int_{0}^{\infty}\tau^{-\frac{x}{x+y}}e^{-\tau}d\tau} \right) g\left( x+y \right) y^{\alpha-1} x^{\beta-1} (x+y)^{\gamma} e^{-x-y}  dy dx=
$$
$$
=\frac{1}{2}\int_{0}^{\infty}\int_{0}^{\infty} f\left( \frac{\int_{0}^{\infty}t^{1+\frac{y}{x+y}-1}e^{-t}dt}{\int_{0}^{\infty}\tau^{1-\frac{x}{x+y}-1}e^{-\tau}d\tau} \right) g\left( x+y \right) y^{\alpha-1} x^{\beta-1} (x+y)^{\gamma} e^{-x-y}  dy dx=
$$
$$
=\frac{1}{2}\int_{0}^{\infty}\int_{0}^{\infty} f\left( \frac{\Gamma\left(1+\frac{y}{x+y}\right)}{\Gamma\left(1-\frac{x}{x+y}\right)} \right) g\left( x+y \right) y^{\alpha-1} x^{\beta-1} (x+y)^{\gamma} e^{-x-y}  dy dx=
$$
$$
=\frac{1}{2}\int_{0}^{\infty}\int_{0}^{\infty} f\left( \frac{\Gamma\left(1+\frac{y}{x+y}\right)}{\Gamma\left(\frac{y}{x+y}\right)} \right) g\left( x+y \right) y^{\alpha-1} x^{\beta-1} (x+y)^{\gamma} e^{-x-y}  dy dx=
$$
$$
=\frac{1}{2}\int_{0}^{\infty}\int_{0}^{\infty} f\left( \frac{y}{x+y} \right) g\left( x+y \right) y^{\alpha-1} x^{\beta-1} (x+y)^{\gamma} e^{-x-y}  dy dx.
$$
Since $\int_{0}^{\infty}g\left( r^{2} \right)r^{2\alpha+2\beta+2\gamma-1}e^{-r^{2}}dr=\frac{1}{2}\Gamma_{g(\cdot)}(\alpha+\beta+\gamma)$, we obtain
$$
\frac{1}{2}\Gamma_{g(\cdot)}(\alpha+\beta+\gamma)B_{f(\cdot)}(\alpha,\beta)=\frac{1}{2}\int_{0}^{\infty}\int_{0}^{\infty} f\left( \frac{y}{x+y} \right) g\left( x+y \right) y^{\alpha-1} x^{\beta-1} (x+y)^{\gamma} e^{-x-y}  dy dx.
$$
$\Box$

Next, we will obtain a similar theorem for the case of an alternating function $Q$. For this we will consider the integrals
\begin{gather}\label{eq1f1}
\Gamma_{|g(\cdot)|}(\alpha+\beta+\gamma)=\int_{0}^{\infty}|g(x)|x^{\alpha+\beta+\gamma-1}e^{-x}dx, \quad \alpha>0, \quad \beta>0, \quad \gamma\geq 0,
\end{gather}
\begin{gather}\label{eq2f1}
B_{|f(\cdot)|}(\alpha,\beta)=\int_{0}^{1}|f(x)|x^{\alpha-1}(1-x)^{\beta-1}dx, \quad \alpha>0, \quad \beta>0,
\end{gather}
$$
I_{\mathbb{R}^{2}_{+}; \left | f\left(\frac{\ast_{_{\alpha}}}{\cdot_{_{\beta}}+\ast_{_{\alpha}}}\right)g\left( \cdot_{_{\beta}}+\ast_{_{\alpha}} \right)\right |}(\alpha,\beta;\gamma):=
$$
\begin{gather}\label{eq3d1f1}
:= \overset {} {\underset {\{T=(x,y)\in \mathbb{R}^{2} \hspace{1mm}  : \hspace{1mm} x\geq 0, y\geq 0 \} }{\int\int} }  \left | f\left(\frac{y}{x+y}\right)g\left( x+y \right)\right | y^{\alpha-1}x^{\beta-1}\left( x+y \right)^{\gamma}e^{-x-y}dydx,
\end{gather}
$$
I^{[1]}_{\mathbb{R}^{2}_{+}; \left|f\left(\frac{\ast_{_{\alpha}}}{\cdot_{_{\beta}}+\ast_{_{\alpha}}}\right)g\left( \cdot_{_{\beta}}+\ast_{_{\alpha}} \right)\right|}(\alpha,\beta;\gamma):=
$$
\begin{gather}\label{eq3f1}
:=\int_{0}^{\infty}\int_{0}^{\infty}\left | f\left(\frac{y}{x+y}\right)g\left( x+y \right)\right |y^{\alpha-1}x^{\beta-1}\left( x+y \right)^{\gamma}e^{-x-y}dydx,
\end{gather}
and
$$
I^{[2]}_{\mathbb{R}^{2}_{+}; \left|f\left(\frac{\ast_{_{\alpha}}}{\cdot_{_{\beta}}+\ast_{_{\alpha}}}\right)g\left( \cdot_{_{\beta}}+\ast_{_{\alpha}} \right)\right|}(\alpha,\beta;\gamma):=
$$
\begin{gather}\label{eq4f1}
:=\int_{0}^{\infty}\int_{0}^{\infty}\left | f\left(\frac{y}{x+y}\right)g\left( x+y \right)\right |y^{\alpha-1}x^{\beta-1}\left( x+y \right)^{\gamma}e^{-x-y}dxdy,
\end{gather}
where $\alpha>0$, $\beta>0$, $\gamma\geq 0$.

{\it Theorem 2.} {\it Let us assume that $Q(x,y)$ is a given function satisfying the following conditions:}

{\it 1) the function} $Q(x,y)$ {\it exists and is continuous at all points in the region } $\{T=(x,y)\in \mathbb{R}^{2} \hspace{1mm}  : \hspace{1mm} x\geq 0, y\geq 0 \} $, {\it with the possible exception of a finite number of points } $M_{1}$, $M_{2}$, ... , $M_{\kappa}$, $\kappa \in \mathbb{N}$, {\it in this region. At points} $M_{1}$, $M_{2}$, ... , $M_{\kappa}$ {\it the function} $Q(x,y)$ {\it tends to} $+ \infty$ or $- \infty$;

{\it 2) integrals} \eqref{eq1f1}, \eqref{eq2f1}, \eqref{eq3d1f1}, \eqref{eq3f1}, \eqref{eq4f1} {\it exist and converge for all values} $\alpha>0$, $\beta>0$, $\gamma\geq 0$;

{\it 3) the equality \eqref{eq5} holds}.

{\it Then the formula \eqref{eq6} will be valid.}

{\it Remark 1.} Theorems 1 and 2 will be valid if in condition 1) we replace the continuity of function $Q(x,y)$ in the area $\{T=(x,y)\in \mathbb{R}^{2} \hspace{1mm}  : \hspace{1mm} x\geq 0, y\geq 0 \} $, with the exception of a finite number of points $M_{1}$, $M_{2}$, ... , $M_{\kappa}$, with continuity of function $Q(x,y)$ in $\{T=(x,y)\in \mathbb{R}^{2} \hspace{1mm}  : \hspace{1mm} x\geq 0, y\geq 0 \} $, with the exception of a finite number of piecewise smooth curves.

Next, we will consider the case of functions that are Riemann integrable.

{\it Theorem 3. Let us assume that $Q(x,y)=f\left(\frac{y}{x+y}\right)g\left( x+y \right)$  is a given function satisfying the following conditions: the functions}
\begin{gather}\label{eq9}
H_{1}(x,y)_{(\alpha,\beta;\gamma)}:= f\left(\frac{y}{x+y}\right)g\left( x+y \right)y^{\alpha-1}x^{\beta-1}\left( x+y \right)^{\gamma}e^{-x-y}, \quad  \alpha>0, \beta>0, \gamma\geq 0,
\end{gather}
\begin{gather}\label{eq10}
H_{2}(x,y)_{(\alpha,\beta;\gamma)}:= f\left(\frac{x}{x+y}\right)g\left( \frac{1}{x}+\frac{1}{y} \right)y^{-\alpha-1}x^{-\beta-1}\left( \frac{1}{x}+\frac{1}{y} \right)^{\gamma}e^{-\frac{1}{x}-\frac{1}{y}},\hspace{1mm} \alpha>0, \beta>0, \gamma\geq 0
\end{gather}
{\it be Riemann  integrable functions in the region} $\{T_{1}=(x,y)\in \mathbb{R}^{2} \hspace{1mm}  : \hspace{1mm} x\geq 0, y\geq 0, x\leq 1, y\leq 1 \}$.

{\it Then the formula \eqref{eq6} will be valid. }

{\it Proof.} The necessity of the conditions of the theorem here is obvious.
$$
I^{[1]}_{\mathbb{R}^{2}_{+}; f\left(\frac{\ast_{_{\alpha}}}{\cdot_{_{\beta}}+\ast_{_{\alpha}}}\right)g\left( \cdot_{_{\beta}}+\ast_{_{\alpha}} \right)}(\alpha,\beta;\gamma)=\int_{0}^{\infty}\int_{0}^{\infty}f\left(\frac{y}{x+y}\right)g\left( x+y \right)y^{\alpha-1}x^{\beta-1}\left( x+y \right)^{\gamma}e^{-x-y}dydx=
$$
$$
=\int_{0}^{1}\int_{0}^{1}f\left(\frac{y}{x+y}\right)g\left( x+y \right)y^{\alpha-1}x^{\beta-1}\left( x+y \right)^{\gamma}e^{-x-y}dydx+
$$
$$
+\int_{0}^{1}\int_{0}^{1}f\left(\frac{x}{x+y}\right)g\left( \frac{1}{x}+\frac{1}{y} \right)y^{-\alpha-1}x^{-\beta-1}\left( \frac{1}{x}+\frac{1}{y} \right)^{\gamma}e^{-\frac{1}{x}-\frac{1}{y}}dydx.
$$
Here we should show that double integrals, when changing the order of integration, will give the same result.
$$
I^{[1]}_{\mathbb{R}^{2}_{+}; f\left(\frac{\ast_{_{\alpha}}}{\cdot_{_{\beta}}+\ast_{_{\alpha}}}\right)g\left( \cdot_{_{\beta}}+\ast_{_{\alpha}} \right)}(\alpha,\beta;\gamma)=\int_{0}^{\infty}\int_{0}^{\infty}f\left(\frac{y}{x+y}\right)g\left( x+y \right)y^{\alpha-1}x^{\beta-1}\left( x+y \right)^{\gamma}e^{-x-y}dydx=
$$
$$
=\int_{0}^{\infty}\int_{0}^{\infty}f\left(\frac{y}{1+y}\right)g\left( x(1+y) \right)y^{\alpha-1}x^{\alpha+\beta-1}\left( x(1+y) \right)^{\gamma}e^{-x(1+y)}dydx=
$$
$$
=\int_{0}^{\infty}\int_{0}^{\infty}f\left(\frac{y}{1+y}\right)g\left( x(1+y) \right)y^{\alpha-1}x^{\alpha+\beta-1}\left( x(1+y) \right)^{\gamma}e^{-x(1+y)}dxdy=
$$
$$
=\int_{0}^{\infty}f\left(\frac{y}{1+y}\right)\frac{y^{\alpha-1}}{(1+y)^{\alpha+\beta}}dy \int_{0}^{\infty}g\left( x \right)x^{\alpha+\beta+\gamma-1}e^{-x}dx=
$$
\begin{gather}\label{eq10f1}
=\Gamma_{g(\cdot)}(\alpha+\beta+\gamma)\int_{0}^{\infty}f\left(\frac{y}{1+y}\right)\frac{y^{\alpha-1}}{(1+y)^{\alpha+\beta}}dy.
\end{gather}
Next, we get
$$
\int_{0}^{\infty}f\left(\frac{y}{1+y}\right)\frac{y^{\alpha-1}}{(1+y)^{\alpha+\beta}}dy=\int_{0}^{\infty}f\left(\frac{1}{1+y}\right)\frac{y^{\beta-1}}{(1+y)^{\alpha+\beta}}dy=
$$
\begin{gather}\label{eq10f2}
=\int_{1}^{\infty}f\left(\frac{1}{y}\right)\frac{(y-1)^{\beta-1}}{y^{\alpha+\beta}}dy=\int_{0}^{1}f\left( y \right)y^{\alpha-1}(1-y)^{\beta-1}dy=B_{f(\cdot)}(\alpha,\beta).
\end{gather}
Substitution \eqref{eq10f2} into \eqref{eq10f1} we obtain formula \eqref{eq6}.

Next, changing the order of integration we obtain
$$
I^{[2]}_{\mathbb{R}^{2}_{+}; f\left(\frac{\ast_{_{\alpha}}}{\cdot_{_{\beta}}+\ast_{_{\alpha}}}\right)g\left( \cdot_{_{\beta}}+\ast_{_{\alpha}} \right)}(\alpha,\beta;\gamma)=\int_{0}^{\infty}\int_{0}^{\infty}f\left(\frac{y}{x+y}\right)g\left( x+y \right)y^{\alpha-1}x^{\beta-1}\left( x+y \right)^{\gamma}e^{-x-y}dxdy=
$$
$$
\int_{0}^{\infty}\int_{0}^{\infty}f\left(\frac{1}{1+x}\right)g\left( y(1+x) \right)x^{\beta-1}y^{\alpha+\beta-1}\left( y(1+x) \right)^{\gamma}e^{-y(1+x)}dxdy=
$$
$$
\int_{0}^{\infty}\int_{0}^{\infty}f\left(\frac{1}{1+x}\right)g\left( y(1+x) \right)x^{\beta-1}y^{\alpha+\beta-1}\left( y(1+x) \right)^{\gamma}e^{-y(1+x)}dydx=
$$
$$
\int_{0}^{\infty}f\left(\frac{1}{1+x}\right)\frac{x^{\beta-1}}{(1+x)^{\alpha+\beta}}\int_{0}^{\infty}g\left( y \right)y^{\alpha+\beta+\gamma-1}e^{-y}dydx=
$$
\begin{gather}\label{eq10f3}
\Gamma_{g(\cdot)}(\alpha+\beta+\gamma)\int_{0}^{\infty}f\left(\frac{1}{1+x}\right)\frac{x^{\beta-1}}{(1+x)^{\alpha+\beta}}dx.
\end{gather}
Next, we get
$$
\int_{0}^{\infty}f\left(\frac{1}{1+x}\right)\frac{x^{\beta-1}}{(1+x)^{\alpha+\beta}}dx=
$$
\begin{gather}\label{eq10f4}
=\int_{1}^{\infty}f\left(\frac{1}{x}\right)\frac{(x-1)^{\beta-1}}{x^{\alpha+\beta}}dx=\int_{0}^{1}f\left( x \right)x^{\alpha-1}(1-x)^{\beta-1}dy=B_{f(\cdot)}(\alpha,\beta).
\end{gather}
Substitution \eqref{eq10f4} into \eqref{eq10f3} we to obtain formula \eqref{eq5}. Thus, under the assumptions of this theorem, an integral of the form \eqref{eq3} does not depend on the order of integration. $\Box$

Next, we will consider the case of a continuous function $Q$.

{\it Theorem 4.} {\it Let us assume that $Q(x,y)=f\left(\frac{y}{x+y}\right)g\left( x+y \right)$  is a given function satisfying the following conditions:}

{\it 1) the finction} $Q(x,y):=f\left(\frac{y}{x+y}\right)g\left( x+y \right)$ {\it exists, is non-negative and is continuous at all points in the region } $\{T=(x,y)\in \mathbb{R}^{2} \hspace{1mm}  : \hspace{1mm} x\geq 0, y\geq 0 \} $;

{\it 2) integrals} \eqref{eq1}, \eqref{eq2}, \eqref{eq3d1}, \eqref{eq3}, \eqref{eq4} {\it exist and converge for all values} $\alpha>0$, $\beta>0$, $\gamma\geq 0$;

{\it 3) the following equality \eqref{eq5} holds};

{\it Then the formula \eqref{eq6} will be valid. }

{\it Proof.} Obviously, Theorem 4 is a special case of Theorem 1. $\Box$

{\it Remark 2.} Theorems 1, 2, 4 and Remark 1 will be valid if in condition 2) of Theorems 1, 2, 4 we require the existence of only integrals \eqref{eq1}, \eqref{eq2} or the existence only integrals  \eqref{eq3d1}, \eqref{eq3}, \eqref{eq4}.

Next, we will introduce the definition of the generalized two-dimensional gamma function in the form of a double integral. One of the main properties of this function is a generalization of the well-known formula for the relationship between the classical gamma function and the beta function.

{\it Definition 1. Let's the function} $Q(x,y)=f\left(\frac{y}{x+y}\right)g\left( x+y \right)$ {\it satisfy the conditions of at least one of Theorems 1-4.
We define the two-dimensional generalized gamma function as a double integral as follows:}
$$
\Gamma_{\mathbb{R}^{2}_{+}; f\left(\frac{\ast_{_{\alpha}}}{\cdot_{_{\beta}}+\ast_{_{\alpha}}}\right)g\left( \cdot_{_{\beta}}+\ast_{_{\alpha}} \right)}(\alpha,\beta;\gamma):=
$$
\begin{gather}\label{eq11}
:=\int_{0}^{\infty}\int_{0}^{\infty}f\left(\frac{y}{x+y}\right)g\left( x+y \right)y^{\alpha-1}x^{\beta-1}\left( x+y \right)^{\gamma}e^{-x-y}dydx, \quad \alpha>0,\hspace{1mm} \beta>0,\hspace{1mm} \gamma\geq 0.
\end{gather}

{\it Theorem 5. (Properties of the two-dimensional generalized gamma function).}

{\it We assume that the conditions of at least one of Theorems 1-4 are satisfied. Then for values of} $\alpha>0$, $\beta>0$, $\gamma\geq 0$ {\it we obtain the following properties of the two-dimensional generalized gamma function:}

1.  $\Gamma_{\mathbb{R}^{2}_{+}; f\left(\frac{\ast_{_{\alpha}}}{\cdot_{_{\beta}}+\ast_{_{\alpha}}}\right)g\left( \cdot_{_{\beta}}+\ast_{_{\alpha}} \right)}(\alpha,\beta;\gamma)=\Gamma_{\mathbb{R}^{2}_{+}; f\left(1-\frac{\cdot_{_{\beta}}}{\cdot_{_{\beta}}+\ast_{_{\alpha}}}\right)g\left( \cdot_{_{\beta}}+\ast_{_{\alpha}} \right)}(\alpha,\beta;\gamma)$;

2.  $\Gamma_{\mathbb{R}^{2}_{+}; f\left(\frac{\ast_{_{\alpha}}}{\cdot_{_{\beta}}+\ast_{_{\alpha}}}\right)g\left( \cdot_{_{\beta}}+\ast_{_{\alpha}} \right)}(\alpha,\beta;\gamma)=\int_{0}^{\infty}\int_{0}^{\infty}f\left(\frac{y}{x+y}\right)g\left( x+y \right)y^{\alpha-1}x^{\beta-1}\left( x+y \right)^{\gamma}e^{-x-y}dxdy=$

$=\overset {} {\underset {\{T=(x,y)\in \mathbb{R}^{2} \hspace{1mm}  : \hspace{1mm} x\geq 0, y\geq 0 \} }{\int\int} } f\left(\frac{y}{x+y}\right)g\left( x+y \right)y^{\alpha-1}x^{\beta-1}\left( x+y \right)^{\gamma}e^{-x-y}dydx$;

3.  $\Gamma_{\mathbb{R}^{2}_{+}; f\left(\frac{\ast_{_{\alpha}}}{\cdot_{_{\beta}}+\ast_{_{\alpha}}}\right)g\left( \cdot_{_{\beta}}+\ast_{_{\alpha}} \right)}(\alpha,\beta;\gamma)=\int_{0}^{\infty}\int_{0}^{\infty}f\left(\frac{x}{x+y}\right)g\left( x+y \right)x^{\alpha-1}y^{\beta-1}\left( x+y \right)^{\gamma}e^{-x-y}dxdy$;

4.  $\Gamma_{\mathbb{R}^{2}_{+}; f\left(\frac{\ast_{_{\alpha}}}{\cdot_{_{\beta}}+\ast_{_{\alpha}}}\right)g\left( \cdot_{_{\beta}}+\ast_{_{\alpha}} \right)}(\alpha,\beta;\gamma)=\int_{0}^{\frac{\pi}{2}}\int_{0}^{\infty}f\left(\frac{\sin (\varphi)}{\sin (\varphi)+\cos (\varphi)}\right)g\left( r\sin (\varphi)+r\cos (\varphi) \right)\cdot$
$$
\cdot (\sin \varphi)^{\alpha-1}(\cos \varphi)^{\beta-1}\left( \sin (\varphi)+\cos (\varphi) \right)^{\gamma}e^{-r(\sin (\varphi)+\cos (\varphi))}r^{\alpha+\beta+\gamma-1} dr d\varphi;
$$

4.1.  $\Gamma_{\mathbb{R}^{2}_{+}; f\left(\frac{\ast_{_{\alpha}}}{\cdot_{_{\beta}}+\ast_{_{\alpha}}}\right)g\left( \cdot_{_{\beta}}+\ast_{_{\alpha}} \right)}(\alpha,\beta;\gamma)=\int_{0}^{\frac{\pi}{2}}\int_{0}^{\infty}f\left(\frac{\cos (\varphi)}{\sin (\varphi)+\cos (\varphi)}\right)g\left( r\sin (\varphi)+r\cos (\varphi) \right)\cdot$
$$
\cdot (\cos \varphi)^{\alpha-1}(\sin \varphi)^{\beta-1}\left( \sin (\varphi)+\cos (\varphi) \right)^{\gamma}e^{-r(\sin (\varphi)+\cos (\varphi))}r^{\alpha+\beta+\gamma-1} dr d\varphi;
$$

5. {\it Generalized formula for the relationship between the gamma function and the beta function}
\begin{gather}\label{eq12}
B_{f(\cdot)}(\alpha,\beta)\Gamma_{g(\cdot)}(\alpha+\beta+\gamma)=\Gamma_{\mathbb{R}^{2}_{+}; f\left(\frac{\ast_{_{\alpha}}}{\cdot_{_{\beta}}+\ast_{_{\alpha}}}\right)g\left( \cdot_{_{\beta}}+\ast_{_{\alpha}} \right)}(\alpha,\beta;\gamma);
\end{gather}

5.1. {\it If we assume} $f\equiv 1$, $g\equiv 1$, $\gamma=0$ {\it in the formula \eqref{eq12}, then we obtain formula for the relationship between the classical gamma function and the beta function}
\begin{gather}\label{eq12a1}
B(\alpha,\beta)\Gamma(\alpha+\beta)=\Gamma(\alpha)\Gamma(\beta);
\end{gather}

5.2. {\it If we assume} $f\equiv 1$, $g\equiv 1$, $\alpha=1$, $\beta=1$ {\it in the formula \eqref{eq12}, then we obtain classical gamma function}
\begin{gather}\label{eq12a1ddgfr}
\Gamma_{\mathbb{R}^{2}_{+};}(1,1;\gamma)=\Gamma(2+\gamma);
\end{gather}

5.3. {\it If we assume} $f\equiv 1$, $g\equiv 1$, $\alpha=\frac{1}{2}$, $\beta=\frac{1}{2}$ {\it in the formula \eqref{eq12}, then we also find a special case of transforming the two-dimensional generalized gamma function into the classical gamma function}
\begin{gather}\label{eq12a1ghlfr}
\Gamma_{\mathbb{R}^{2}_{+};}(1/2,1/2;\gamma)=\pi\Gamma(1+\gamma);
\end{gather}

5.4. {\it If we assume} $f\equiv 1$, $g\equiv 1$ {\it in the formula \eqref{eq12}, then we obtain }
\begin{gather}\label{eq12a34rt}
\int_{0}^{\infty}\int_{0}^{\infty}y^{\alpha-1}x^{\beta-1}\left( x+y \right)^{\gamma}e^{-x-y}dydx=B(\alpha,\beta)\Gamma(\alpha+\beta+\gamma);
\end{gather}

6. {\it If the function} $g$ {\it identically not equal to zero, then}
$$
\frac{\Gamma_{g(\cdot)}(\alpha+\beta+\gamma)}{\Gamma_{g(\cdot)}(\alpha+\beta+\gamma+1)}=\frac{\Gamma_{\mathbb{R}^{2}_{+}; f\left(\frac{\ast_{_{\alpha}}}{\cdot_{_{\beta}}+\ast_{_{\alpha}}}\right)g\left( \cdot_{_{\beta}}+\ast_{_{\alpha}} \right)}(\alpha,\beta;\gamma)}{\Gamma_{\mathbb{R}^{2}_{+}; f\left(\frac{\ast_{_{\alpha}}}{\cdot_{_{\beta}}+\ast_{_{\alpha}}}\right)g\left( \cdot_{_{\beta}}+\ast_{_{\alpha}} \right)}(\alpha,\beta;\gamma+1)}.
$$

7. {\it If the limit is valid} $\lim e^{-x}x^{\alpha+\beta+\gamma-1}g(x)=0$, {\it if} $x\rightarrow 0$ {\it and} $x\rightarrow +\infty$, {\it then we find}
$$
\Gamma_{\mathbb{R}^{2}_{+}; f\left(\frac{\ast_{_{\alpha}}}{\cdot_{_{\beta}}+\ast_{_{\alpha}}}\right)g\left( \cdot_{_{\beta}}+\ast_{_{\alpha}} \right)}(\alpha,\beta;\gamma+1)=(\alpha+\beta+\gamma)\Gamma_{\mathbb{R}^{2}_{+}; f\left(\frac{\ast_{_{\alpha}}}{\cdot_{_{\beta}}+\ast_{_{\alpha}}}\right)g\left( \cdot_{_{\beta}}+\ast_{_{\alpha}} \right)}(\alpha,\beta;\gamma)+
$$
$$
+\Gamma_{\mathbb{R}^{2}_{+}; f\left(\frac{\ast_{_{\alpha}}}{\cdot_{_{\beta}}+\ast_{_{\alpha}}}\right)g' \left( \cdot_{_{\beta}}+\ast_{_{\alpha}} \right)}(\alpha,\beta;\gamma+1).
$$

{\it Proof.} Properties 1-4 are obvious. Property 5 follows from {\it Definition 1} of the two-dimensional generalized gamma function \eqref{eq11} and the formula \eqref{eq6}, due to the fulfillment of the conditions of at least one of {\it Theorems 1-4.} Property 7 follows from Property 5 and formula
$$
\Gamma_{ g (\cdot)}(\alpha+\beta+\gamma+1)=(\alpha+\beta+\gamma)\Gamma_{ g (\cdot)}(\alpha+\beta+\gamma)+\Gamma_{g' ( \cdot)}(\alpha+\beta+\gamma+1),\quad \alpha>0,  \hspace{1mm} \beta>0, \hspace{1mm} \gamma\geq 0.
$$

$\Box$

{\it Lemma 1. Let us assume that following conditions are met:

1) the function} $f'(x)$ {\it exists on} $(0,1)$ {\it and the following integral exists}
\begin{gather}\label{eq20}
\int_{0}^{1}f'(x)x^{\alpha}(1-x)^{\beta}dx, \quad \alpha>0, \quad \beta > 0;
\end{gather}

{\it 2) the following limits are valid}
$$
\overset {} {\underset {x \rightarrow 0+ }{\lim} } f(x)x^{\alpha}(1-x)^{\beta}=0, \quad \alpha>0, \quad \beta>0,
$$
$$
\overset {} {\underset {x \rightarrow 1- }{\lim} } f(x)x^{\alpha}(1-x)^{\beta}=0, \quad \alpha>0, \quad \beta>0.
$$
{\it Then the following formula will be valid}
\begin{gather}\label{eq21}
B_{f(\cdot)}(\alpha,\beta)-B_{f(\cdot)}(\alpha+1,\beta)-B_{f(\cdot)}(\alpha,\beta+1)=0, \quad \alpha>0, \quad \beta>0.
\end{gather}
{\it Proof.} We obtain
\begin{gather}\label{eq22}
\int_{0}^{1}f'(x)x^{\alpha}(1-x)^{\beta}dx=-\alpha B_{f(\cdot)}(\alpha,\beta+1)+\beta B_{f(\cdot)}(\alpha+1,\beta), \quad \alpha>0, \quad \beta > 0.
\end{gather}
Next, we get
$$
\int_{0}^{1}f'(x)x^{\alpha}(1-x)^{\beta}dx=-\int_{0}^{1}f(x)\left( \alpha x^{\alpha-1}(1-x)^{\beta} +\beta x^{\alpha}(1-x)^{\beta-1}\right)dx=
$$
$$
=-\int_{0}^{1}f(x)x^{\alpha-1}(1-x)^{\beta-1}\left( \alpha (1-x) +\beta x \right)dx=
$$
$$
=-\alpha B_{f(\cdot)}(\alpha,\beta)-(\beta-\alpha)B_{f(\cdot)}(\alpha+1,\beta), \quad \alpha>0, \quad \beta > 0.
$$
Then we obtain
\begin{gather}\label{eq23}
\int_{0}^{1}f'(x)x^{\alpha}(1-x)^{\beta}dx=(\alpha+\beta)B_{f(\cdot)}(\alpha+1,\beta)-\alpha B_{f(\cdot)}(\alpha,\beta), \quad \alpha>0, \quad \beta > 0.
\end{gather}
Equating the right sides of formulas \eqref{eq22} end \eqref{eq23}, we get \eqref{eq21}.
$\Box$

{\it Lemma 2. Let us assume that following conditions are met:

1) the function} $f'(x)$ {\it exists on} $(0,1)$ {\it and following integral exists}
$$
\int_{0}^{1}f'(x)x^{\alpha}(1-x)^{\beta}dx, \quad \alpha>0, \quad \beta > 0;
$$

{\it 2) the following limits are valid}
$$
\overset {} {\underset {x \rightarrow 0+ }{\lim} } f(x)x^{\alpha}(1-x)^{\beta}=0, \quad \alpha>0, \quad \beta>0,
$$
$$
\overset {} {\underset {x \rightarrow 1- }{\lim} } f(x)x^{\alpha}(1-x)^{\beta}=0, \quad \alpha>0, \quad \beta>0;
$$

{\it 3) the function} $Q(x,y)=f\left(\frac{y}{x+y}\right)g\left( x+y \right)$ {\it satisfy the conditions of at least one of Theorems 1-4.}

{\it Then the following formulas will be valid}
$$
\Gamma_{\mathbb{R}^{2}_{+}; f\left(\frac{\ast_{_{\alpha}}}{\cdot_{_{\beta}}+\ast_{_{\alpha}}}\right)g\left( \cdot_{_{\beta}}+\ast_{_{\alpha}} \right)}(\alpha+1,\beta;\gamma)+\Gamma_{\mathbb{R}^{2}_{+}; f\left(\frac{\ast_{_{\alpha}}}{\cdot_{_{\beta}}+\ast_{_{\alpha}}}\right)g\left( \cdot_{_{\beta}}+\ast_{_{\alpha}} \right)}(\alpha,\beta+1;\gamma)=
$$
\begin{gather}\label{eq24}
=\Gamma_{\mathbb{R}^{2}_{+}; f\left(\frac{\ast_{_{\alpha}}}{\cdot_{_{\beta}}+\ast_{_{\alpha}}}\right)g\left( \cdot_{_{\beta}}+\ast_{_{\alpha}} \right)}(\alpha,\beta;\gamma+1), \quad \alpha>0, \quad \beta>0, \quad \gamma\geq 0;
\end{gather}
$$
\beta\Gamma_{\mathbb{R}^{2}_{+}; f\left(\frac{\ast_{_{\alpha}}}{\cdot_{_{\beta}}+\ast_{_{\alpha}}}\right)g\left( \cdot_{_{\beta}}+\ast_{_{\alpha}} \right)}(\alpha+1,\beta;\gamma+1)-\alpha \Gamma_{\mathbb{R}^{2}_{+}; f\left(\frac{\ast_{_{\alpha}}}{\cdot_{_{\beta}}+\ast_{_{\alpha}}}\right)g\left( \cdot_{_{\beta}}+\ast_{_{\alpha}} \right)}(\alpha,\beta+1;\gamma+1)=
$$
\begin{gather}\label{eq24sdr455}
= \Gamma_{\mathbb{R}^{2}_{+}; f'\left(\frac{\ast_{_{\alpha}}}{\cdot_{_{\beta}}+\ast_{_{\alpha}}}\right)g\left( \cdot_{_{\beta}}+\ast_{_{\alpha}} \right)}(\alpha+1,\beta+1;\gamma), \quad \alpha>0, \quad \beta>0, \quad \gamma\geq 0.
\end{gather}

{\it Proof.} From equality \eqref{eq21} we find
$$
B_{f(\cdot)}(\alpha,\beta)\Gamma_{g(\cdot)}(\alpha+\beta+\gamma)=\Gamma_{g(\cdot)}(\alpha+\beta+\gamma)\left( B_{f(\cdot)}(\alpha+1,\beta)+B_{f(\cdot)}(\alpha,\beta+1)\right),
$$
and futher, from equality \eqref{eq11}, we get
$$
\Gamma_{\mathbb{R}^{2}_{+}; f\left(\frac{\ast_{_{\alpha}}}{\cdot_{_{\beta}}+\ast_{_{\alpha}}}\right)g\left( \cdot_{_{\beta}}+\ast_{_{\alpha}} \right)}(\alpha,\beta;\gamma)=\Gamma_{g(\cdot)}(\alpha+\beta+\gamma)\left( B_{f(\cdot)}(\alpha+1,\beta)+B_{f(\cdot)}(\alpha,\beta+1)\right)=
$$
$$
=\frac{\Gamma_{g(\cdot)}(\alpha+\beta+\gamma)}{\Gamma_{g(\cdot)}(\alpha+\beta+\gamma+1)}\left( \Gamma_{\mathbb{R}^{2}_{+}; f\left(\frac{\ast_{_{\alpha}}}{\cdot_{_{\beta}}+\ast_{_{\alpha}}}\right)g\left( \cdot_{_{\beta}}+\ast_{_{\alpha}} \right)}(\alpha+1,\beta;\gamma)+ \Gamma_{\mathbb{R}^{2}_{+}; f\left(\frac{\ast_{_{\alpha}}}{\cdot_{_{\beta}}+\ast_{_{\alpha}}}\right)g\left( \cdot_{_{\beta}}+\ast_{_{\alpha}} \right)}(\alpha,\beta+1;\gamma) \right).
$$
From the last equality we find
$$
\frac{\Gamma_{\mathbb{R}^{2}_{+}; f\left(\frac{\ast_{_{\alpha}}}{\cdot_{_{\beta}}+\ast_{_{\alpha}}}\right)g\left( \cdot_{_{\beta}}+\ast_{_{\alpha}} \right)}(\alpha,\beta;\gamma)}{\Gamma_{\mathbb{R}^{2}_{+}; f\left(\frac{\ast_{_{\alpha}}}{\cdot_{_{\beta}}+\ast_{_{\alpha}}}\right)g\left( \cdot_{_{\beta}}+\ast_{_{\alpha}} \right)}(\alpha+1,\beta;\gamma)+ \Gamma_{\mathbb{R}^{2}_{+}; f\left(\frac{\ast_{_{\alpha}}}{\cdot_{_{\beta}}+\ast_{_{\alpha}}}\right)g\left( \cdot_{_{\beta}}+\ast_{_{\alpha}} \right)}(\alpha,\beta+1;\gamma)}=
$$
$$
=\frac{\Gamma_{g(\cdot)}(\alpha+\beta+\gamma)B_{f(\cdot)}(\alpha,\beta)}{\Gamma_{g(\cdot)}(\alpha+\beta+\gamma+1)B_{f(\cdot)}(\alpha,\beta)}=\frac{\Gamma_{\mathbb{R}^{2}_{+}; f\left(\frac{\ast_{_{\alpha}}}{\cdot_{_{\beta}}+\ast_{_{\alpha}}}\right)g\left( \cdot_{_{\beta}}+\ast_{_{\alpha}} \right)}(\alpha,\beta;\gamma)}{\Gamma_{\mathbb{R}^{2}_{+}; f\left(\frac{\ast_{_{\alpha}}}{\cdot_{_{\beta}}+\ast_{_{\alpha}}}\right)g\left( \cdot_{_{\beta}}+\ast_{_{\alpha}} \right)}(\alpha,\beta;\gamma+1)},
$$
from which we get formula \eqref{eq24}. Formula \eqref{eq24sdr455} we obtain from the formula \eqref{eq22}.
$\Box$

Next, we will introduce differentiation formulas for our generalization of the gamma function.

{\it Lemma 3. Let's the function} $Q(x,y)=f\left(\frac{y}{x+y}\right)g\left( x+y \right)$ {\it satisfy the conditions of at least one of Theorems 1-4.

Then for} $\forall$ $\alpha>0$, $\beta>0$, $\gamma\geq 0$ {\it and for} $\forall$ $l,m,n$, {\it where} $\{l,m,n\}\in \mathbb{N}\cup \{0\}$,  {\it exist derivatives of all orders for the two-dimensional generalized gamma function } $\Gamma_{\mathbb{R}^{2}_{+}; f\left(\frac{\ast_{_{\alpha}}}{\cdot_{_{\beta}}+\ast_{_{\alpha}}}\right)g\left( \cdot_{_{\beta}}+\ast_{_{\alpha}} \right)}(\alpha,\beta;\gamma)$ {\it and the following formula is valid}
$$
\frac{\partial^{l}}{\partial \gamma^{l}}\frac{\partial^{n}}{\partial \beta^{n}}\frac{\partial^{m}}{\partial \alpha^{m}}\Gamma_{\mathbb{R}^{2}_{+}; f\left(\frac{\ast_{_{\alpha}}}{\cdot_{_{\beta}}+\ast_{_{\alpha}}}\right)g\left( \cdot_{_{\beta}}+\ast_{_{\alpha}} \right)}(\alpha,\beta;\gamma)=
$$
$$
=\int_{0}^{\infty}\int_{0}^{\infty}f\left(\frac{y}{x+y}\right)g\left( x+y \right)\left( \log(x+y) \right)^{l}(\log (x))^{n}(\log (y))^{m}y^{\alpha-1}x^{\beta-1}\left( x+y \right)^{\gamma}e^{-x-y}dydx,
$$
{\it where } $\{l,m,n\}\subset \mathbb{N}\cup \{0\}$, $\alpha>0$, $\beta>0$, $\gamma\geq 0$.

{\it Theorem 6. (Formula for differentiating the two-dimensional generalized gamma function)}

{\it Let's the function} $Q(x,y)=f\left(\frac{y}{x+y}\right)g\left( x+y \right)$ {\it satisfy the conditions of at least one of Theorems 1-4. Then the following formula will be valid}
$$
\int_{0}^{\infty}\int_{0}^{\infty}f\left(\frac{y}{x+y}\right)g\left( x+y \right)\left( \log(x+y) \right)^{l}(\log (x))^{n}(\log (y))^{m}y^{\alpha-1}x^{\beta-1}\left( x+y \right)^{\gamma}e^{-x-y}dydx=
$$
$$
= \sum_{i=0}^{n} \sum_{j=0}^{m} {n \choose i} {m \choose j} \int_{0}^{\infty}  g\left( r \right) \left( \log \left(r\right) \right)^{i+j+l}r^{\alpha+\beta+\gamma-1} e^{-r} dr \cdot
$$
\begin{gather}\label{eq26}
\cdot \int_{0}^{1} f\left( x \right) \left(\log\left( x \right)\right)^{m-j}    \left(\log\left( 1-x \right)\right)^{n-i} x^{\alpha-1} (1-x)^{\beta-1} dx,
\end{gather}
{\it where } $\{l,m,n\}\subset \mathbb{N}\cup \{0\}$, $\alpha>0$, $\beta>0$, $\gamma\geq 0$, ${N \choose M}:=\frac{N!}{M!(N-M)!}$. {\it In the notation style,
according to Definition 1, our differentiating formula} \eqref{eq26} {\it has the form}
$$
\frac{\partial^{l}}{\partial \gamma^{l}}\frac{\partial^{n}}{\partial \beta^{n}}\frac{\partial^{m}}{\partial \alpha^{m}}\Gamma_{\mathbb{R}^{2}_{+}; f\left(\frac{\ast_{_{\alpha}}}{\cdot_{_{\beta}}+\ast_{_{\alpha}}}\right)g\left( \cdot_{_{\beta}}+\ast_{_{\alpha}} \right)}(\alpha,\beta;\gamma)=
$$
\begin{gather}\label{eq27}
= \sum_{j=0}^{m}\sum_{i=0}^{n} {m \choose j}{n \choose i}\Gamma_{g(\cdot)}^{(i+j+l)}(\alpha+\beta+\gamma)\frac{\partial^{n-i}}{\partial \beta^{n-i}}\frac{\partial^{m-j}}{\partial \alpha^{m-j}}B_{f(\cdot)}(\alpha,\beta),
\end{gather}
{\it where } $\{l,m,n\}\subset \mathbb{N}\cup \{0\}$, $\alpha>0$, $\beta>0$, $\gamma\geq 0$.

{\it Proof.}
$$
\frac{\partial^{l}}{\partial \gamma^{l}}\frac{\partial^{n}}{\partial \beta^{n}}\frac{\partial^{m}}{\partial \alpha^{m}}\Gamma_{\mathbb{R}^{2}_{+}; f\left(\frac{\ast_{_{\alpha}}}{\cdot_{_{\beta}}+\ast_{_{\alpha}}}\right)g\left( \cdot_{_{\beta}}+\ast_{_{\alpha}} \right)}(\alpha,\beta;\gamma)=
$$
$$
=\int_{0}^{\infty}\int_{0}^{\infty}f\left(\frac{y}{x+y}\right)g\left( x+y \right)\left( \log(x+y) \right)^{l}(\log (x))^{n}(\log (y))^{m}y^{\alpha-1}x^{\beta-1}\left( x+y \right)^{\gamma}e^{-x-y}dydx=
$$
$$
=4\int_{0}^{\infty}\int_{0}^{\infty} f\left(\frac{y^{2}}{x^{2}+y^{2}}\right) g\left( x^{2}+y^{2} \right) \left( \log\left(x^{2}+y^{2}\right) \right)^{l} \cdot
$$
$$
\cdot \left(\log \left(x^{2}\right)\right)^{n} \left(\log \left(y^{2}\right)\right)^{m} y^{2\alpha-1} x^{2\beta-1} \left(x^{2}+y^{2}\right)^{\gamma} e^{-x^{2}-y^{2}}  dy dx =
$$
$$
=4\int_{0}^{\infty}\int_{0}^{\infty} f\left( \frac{(x^{2}+y^{2})\int_{0}^{\infty}t^{y^{2}}e^{-t^{x^{2}+y^{2}}}t^{x^{2}+y^{2}-1}dt}{(x^{2}+y^{2})\int_{0}^{\infty}\tau^{-x^{2}}e^{-\tau^{x^{2}+y^{2}}}\tau^{x^{2}+y^{2}-1}d\tau} \right) g\left( x^{2}+y^{2} \right) \cdot
$$
$$
\cdot \left( \log\left(x^{2}+y^{2}\right) \right)^{l} \left(\log \left(x^{2}\right)\right)^{n} \left(\log \left(y^{2}\right)\right)^{m} y^{2\alpha-1} x^{2\beta-1} \left(x^{2}+y^{2}\right)^{\gamma} e^{-x^{2}-y^{2}}  dy dx =
$$
$$
=4\int_{0}^{\frac{\pi}{2}}\int_{0}^{\infty} f\left( \frac{r^{2}\int_{0}^{\infty}t^{r^{2}\sin^{2}\varphi}e^{-t^{r^{2}}}t^{r^{2}-1}dt}{r^{2}\int_{0}^{\infty}\tau^{-r^{2}\cos^{2}\varphi}e^{-\tau^{r^{2}}}\tau^{r^{2}-1}d\tau} \right) g\left( r^{2} \right)\cdot
$$
$$
\cdot \left( \log \left(r^{2}\right) \right)^{l} \left(\log\left( r^{2} \cos^{2}\varphi \right)\right)^{n} \left(\log\left( r^{2} \sin^{2}\varphi \right)\right)^{m} (r \sin\varphi)^{2\alpha-1} (r \cos\varphi)^{2\beta-1} r^{2\gamma} e^{-r^{2}} r dr d\varphi=
$$
$$
=4\int_{0}^{\frac{\pi}{2}}\int_{0}^{\infty} f\left( \frac{\int_{0}^{\infty}t^{\sin^{2}\varphi}e^{-t}dt}{\int_{0}^{\infty}\tau^{-\cos^{2}\varphi}e^{-\tau}d\tau} \right) g\left( r^{2} \right) (r \sin\varphi)^{2\alpha-1} (r \cos\varphi)^{2\beta-1} r^{2\gamma} \left( \log \left(r^{2}\right) \right)^{l} \cdot
$$
$$
\cdot \sum_{i=0}^{n} {n \choose i} \left(\log\left( \cos^{2}\varphi \right)\right)^{n-i} \left( \log \left(r^{2}\right) \right)^{i} \sum_{j=0}^{m} {m \choose j} \left(\log\left( \sin^{2}\varphi \right)\right)^{m-j}\left( \log \left(r^{2}\right) \right)^{j}  e^{-r^{2}} r dr d\varphi=
$$
$$
=4 \sum_{i=0}^{n} \sum_{j=0}^{m} {n \choose i} {m \choose j} \int_{0}^{\frac{\pi}{2}}\int_{0}^{\infty} f\left( \frac{\int_{0}^{\infty}t^{\sin^{2}\varphi}e^{-t}dt}{\int_{0}^{\infty}\tau^{-\cos^{2}\varphi}e^{-\tau}d\tau} \right) g\left( r^{2} \right) \left( \log \left(r^{2}\right) \right)^{i+j+l} \cdot
$$
$$
\cdot \left(\log\left( \cos^{2}\varphi \right)\right)^{n-i}    \left(\log\left( \sin^{2}\varphi \right)\right)^{m-j} (r \sin\varphi)^{2\alpha-1} (r \cos\varphi)^{2\beta-1} r^{2\gamma} e^{-r^{2}} r dr d\varphi=
$$
$$
=4 \sum_{i=0}^{n} \sum_{j=0}^{m} {n \choose i} {m \choose j} \int_{0}^{\infty}  g\left( r^{2} \right) \left( \log \left(r^{2}\right) \right)^{i+j+l}r^{2(\alpha+\beta+\gamma)-1} e^{-r^{2}} dr  \cdot
$$
$$
\cdot \int_{0}^{\frac{\pi}{2}} f\left( \sin^{2}\varphi \right) \left(\log\left( \cos^{2}\varphi \right)\right)^{n-i}    \left(\log\left( \sin^{2}\varphi \right)\right)^{m-j} (\sin\varphi)^{2\alpha-1} (\cos\varphi)^{2\beta-1} d\varphi=
$$
$$
= \sum_{i=0}^{n} \sum_{j=0}^{m} {n \choose i} {m \choose j} \int_{0}^{\infty}  g\left( r \right) \left( \log \left(r\right) \right)^{i+j+l}r^{\alpha+\beta+\gamma-1} e^{-r} dr  \cdot
$$
$$
\cdot \int_{0}^{1} f\left( x \right) \left(\log\left( x \right)\right)^{m-j}    \left(\log\left( 1-x \right)\right)^{n-i} x^{\alpha-1} (1-x)^{\beta-1} dx=
$$
$$
= \sum_{j=0}^{m}\sum_{i=0}^{n} {m \choose j}{n \choose i}\Gamma_{g(\cdot)}^{(i+j+l)}(\alpha+\beta+\gamma)\frac{\partial^{n-i}}{\partial \beta^{n-i}}\frac{\partial^{m-j}}{\partial \alpha^{m-j}}B_{f(\cdot)}(\alpha,\beta).
$$
In this proof we applied the transformations
$$
\frac{y^{2}}{x^{2}+y^{2}}=\frac{\Gamma\left( \frac{y^{2}}{x^{2}+y^{2}}+1 \right)}{\Gamma\left( \frac{y^{2}}{x^{2}+y^{2}} \right)}=\frac{\int_{0}^{\infty}t^{\frac{y^{2}}{x^{2}+y^{2}}}e^{-t}dt}{\int_{0}^{\infty}\tau^{-\frac{x^{2}}{x^{2}+y^{2}}}e^{-\tau}d\tau}=\frac{(x^{2}+y^{2})\int_{0}^{\infty}t^{y^{2}}e^{-t^{x^{2}+y^{2}}}t^{x^{2}+y^{2}-1}dt}{(x^{2}+y^{2})\int_{0}^{\infty}\tau^{-x^{2}}e^{-\tau^{x^{2}+y^{2}}}\tau^{x^{2}+y^{2}-1}d\tau},
$$
$$
f\left( \frac{\int_{0}^{\infty}t^{\sin^{2}\varphi}e^{-t}dt}{\int_{0}^{\infty}\tau^{-\cos^{2}\varphi}e^{-\tau}d\tau} \right)=f\left( \frac{\Gamma(1+\sin ^{2}\varphi)}{\Gamma(\sin ^{2}\varphi)} \right)=f\left( \sin ^{2}\varphi \right), \quad \varphi \in \left( 0,\frac{\pi}{2}\right).
$$

$\Box$

{\it Remark 3.} Using the equality $B_{f(\cdot)}(\alpha,\beta)=B_{f(1-\cdot)}(\beta,\alpha)$, from formula \eqref{eq27}, we find
$$
\int_{0}^{\infty}\int_{0}^{\infty}f\left(\frac{x}{x+y}\right)g\left( x+y \right)\left( \log(x+y) \right)^{l}(\log (x))^{n}(\log (y))^{m}y^{\alpha-1}x^{\beta-1}\left( x+y \right)^{\gamma}e^{-x-y}dydx=
$$
$$
= \sum_{j=0}^{m}\sum_{i=0}^{n} {m \choose j}{n \choose i}\Gamma_{g(\cdot)}^{(i+j+l)}(\alpha+\beta+\gamma)\frac{\partial^{n-i}}{\partial \beta^{n-i}}\frac{\partial^{m-j}}{\partial \alpha^{m-j}}B_{f(\cdot)}(\beta,\alpha),
$$
where $\{l,m,n\}\subset \mathbb{N}\cup \{0\}$, $\alpha>0$, $\beta>0$, $\gamma\geq 0$.

{\it Remark 4.} Assuming $f\equiv 1$, $g\equiv 1$, $\gamma=0$ in the formula \eqref {eq27}, we obtain a formula for differentiating the classical gamma and beta functions of any order
\begin{gather}\label{eq101k1}
\Gamma^{(n)}(\beta)\Gamma^{(m)}(\alpha)= \sum_{j=0}^{m}\sum_{i=0}^{n} {m \choose j}{n \choose i}\Gamma^{(i+j)}(\alpha+\beta)\frac{\partial^{n-i}}{\partial \beta^{n-i}}\frac{\partial^{m-j}}{\partial \alpha^{m-j}}B(\alpha,\beta),
\end{gather}
where $\{m,n\}\subset \mathbb{N}\cup \{0\}$, $\alpha>0$, $\beta>0$. The formula \eqref{eq101k1} is also a generalization of formula for the relationship between the classical gamma function and the beta function $\Gamma(\alpha)\Gamma(\beta)=\Gamma(\alpha+\beta)B(\alpha,\beta)$, if $m=n=0$.

{\it Remark 5.} We will assume that $\varrho_{1}, \sigma_{1}, \varrho_{2}, \sigma_{2}, \varrho_{3}, \sigma_{3} \in \mathbb{R}$, then there are all these values belonging to the set of real numbers. Theorems 1-6, Lems 1-3 and Remark 1-4 will be valid if we consider our values of $\alpha$, $\beta$, $\gamma$ in the complex plane, provided that $\alpha:=\varrho_{1}+i\sigma_{1}, \varrho_{1}>0$; $\beta:=\varrho_{2}+i\sigma_{2}, \varrho_{2}>0$; $\gamma:=\varrho_{3}+i\sigma_{3}, \varrho_{3}\geq 0$.

\medskip
\section{Practical applications of the generalized two-dimensional gamma function for calculating double integrals, including those containing some special functions.}

In this chapter, we will present several theorems containing integral equalities using the two-dimensional generalized gamma function.

{\it Theorem 7.}
{\it Let's the function} $f\left(\frac{y}{x+y}\right)$ {\it satisfy the conditions of at least one of Theorems 1-4, if} $g\equiv 1$. {\it Then for any} $\alpha>0$, $\beta>0$, $\gamma\geq 0$ {\it the following formulas will be valid:}
$$
\int_{0}^{\infty}\int_{0}^{\infty}f\left(\frac{y}{x+y}\right)y^{\alpha-1}x^{\beta-1}\left( x+y \right)^{\gamma}e^{-a(x+y)}\cos(b(x+y))dydx=
$$
\begin{gather}\label{eq101f}
=\frac{\Gamma(\alpha+\beta+\gamma)}{\left(\sqrt{a^{2}+b^{2}}\right)^{\alpha+\beta+\gamma}}\cos\left((\alpha+\beta+\gamma)\arctan\frac{b}{a}\right)B_{f(\cdot)}(\alpha,\beta), \quad a>0;
\end{gather}
$$
\int_{0}^{\infty}\int_{0}^{\infty}f\left(\frac{y}{x+y}\right)y^{\alpha-1}x^{\beta-1}\left( x+y \right)^{\gamma}e^{-a(x+y)}\sin(b(x+y))dydx=
$$
\begin{gather}\label{eq102f}
=\frac{\Gamma(\alpha+\beta+\gamma)}{\left(\sqrt{a^{2}+b^{2}}\right)^{\alpha+\beta+\gamma}}\sin\left((\alpha+\beta+\gamma)\arctan\frac{b}{a}\right)B_{f(\cdot)}(\alpha,\beta), \quad a>0.
\end{gather}
{\it Proof.} We obtain
$$
B_{f(\cdot)}(\alpha,\beta)\int_{0}^{\infty}r^{2\alpha+2\beta+2\gamma-1}e^{-(a+bi)r^{2}}dr=
$$
$$
=2\int_{0}^{\frac{\pi}{2}}\int_{0}^{\infty} f\left( \frac{\int_{0}^{\infty}t^{\sin^{2}\varphi}e^{-t}dt}{\int_{0}^{\infty}\tau^{-\cos^{2}\varphi}e^{-\tau}d\tau} \right) (r \sin\varphi)^{2\alpha-1} (r \cos\varphi)^{2\beta-1} r^{2\gamma} e^{-(a+bi)r^{2}} r dr d\varphi=
$$
$$
=\frac{1}{2}\int_{0}^{\infty}\int_{0}^{\infty} f\left( \frac{y}{x+y} \right) y^{\alpha-1} x^{\beta-1} (x+y)^{\gamma} e^{-(a+bi)(x+y)}  dy dx=
$$
$$
=\frac{1}{2}\int_{0}^{\infty}\int_{0}^{\infty} f\left( \frac{y}{x+y} \right) y^{\alpha-1} x^{\beta-1} (x+y)^{\gamma} e^{-a(x+y)}\left(\cos(b(x+y))-i\sin(b(x+y))\right)  dy dx=
$$
$$
=\frac{1}{2}\int_{0}^{\infty}\int_{0}^{\infty} f\left( \frac{y}{x+y} \right) y^{\alpha-1} x^{\beta-1} (x+y)^{\gamma} e^{-a(x+y)}\cos(b(x+y))  dy dx-
$$
\begin{gather}\label{eq103f}
-i\frac{1}{2}\int_{0}^{\infty}\int_{0}^{\infty} f\left( \frac{y}{x+y} \right) y^{\alpha-1} x^{\beta-1} (x+y)^{\gamma} e^{-a(x+y)}\sin(b(x+y))  dy dx.
\end{gather}

And further, we find
$$
B_{f(\cdot)}(\alpha,\beta)\int_{0}^{\infty}r^{2\alpha+2\beta+2\gamma-1}e^{-(a+bi)r^{2}}dr=
$$
$$
=B_{f(\cdot)}(\alpha,\beta)\int_{0}^{\infty}r^{2\alpha+2\beta+2\gamma-1}e^{-ar^{2}}\cos(br)dr-i B_{f(\cdot)}(\alpha,\beta)\int_{0}^{\infty}r^{2\alpha+2\beta+2\gamma-1}e^{-ar^{2}}\sin(br)dr=
$$
\begin{gather}\label{eq104f}
=\frac{1}{2}B_{f(\cdot)}(\alpha,\beta)\int_{0}^{\infty}r^{\alpha+\beta+\gamma-1}e^{-ar}\cos(br)dr-i\frac{1}{2}B_{f(\cdot)}(\alpha,\beta)\int_{0}^{\infty}r^{\alpha+\beta+\gamma-1}e^{-ar}\sin(br)dr.
\end{gather}
Equating the real and imaginary parts of equalities \eqref{eq103f}, \eqref{eq104f}, we obtain \eqref{eq101f}, \eqref{eq102f}, with further consideration of the following two equalities
$$
\int_{0}^{\infty}r^{\alpha+\beta+\gamma-1}e^{-ar}\sin(br)dr=\frac{\Gamma(\alpha+\beta+\gamma)}{\left(\sqrt{a^{2}+b^{2}}\right)^{\alpha+\beta+\gamma}}\sin\left((\alpha+\beta+\gamma)\arctan\frac{b}{a}\right), \quad a>0;
$$
$$
\int_{0}^{\infty}r^{\alpha+\beta+\gamma-1}e^{-ar}\cos(br)dr=\frac{\Gamma(\alpha+\beta+\gamma)}{\left(\sqrt{a^{2}+b^{2}}\right)^{\alpha+\beta+\gamma}}\cos\left((\alpha+\beta+\gamma)\arctan\frac{b}{a}\right), \quad a>0.
$$
$\Box$

{\it Theorem 8.}
{\it Let's the function} $Q_{1}(x,y):=f\left(\frac{y}{x+y}\right)g\left( x+y \right)e^{(1-a)(x+y)}$ {\it satisfy the conditions of function} $Q(x,y)$ {\it of at least one of Theorems 1-4. Then for any} $\alpha>0$, $\beta>0$, $\gamma\geq 0$, $a>0$ {\it the following formulas will be valid:}
$$
\int_{0}^{\infty}\int_{0}^{\infty}f\left(\frac{y}{x+y}\right)g\left( x+y \right)y^{\alpha-1}x^{\beta-1}\left( x+y \right)^{\gamma}e^{-a(x+y)}\cos(b(x+y))dydx=
$$
$$
=B_{f(\cdot)}(\alpha,\beta)\int_{0}^{\infty}g(r)r^{\alpha+\beta+\gamma-1}e^{-ar}\cos(br)dr=B_{f(\cdot)}(\alpha,\beta)\Gamma_{ \cos(b\cdot)e^{\cdot(1-a)}g(\cdot) }(\alpha+\beta+\gamma);
$$
$$
\int_{0}^{\infty}\int_{0}^{\infty}f\left(\frac{y}{x+y}\right)g\left( x+y \right)y^{\alpha-1}x^{\beta-1}\left( x+y \right)^{\gamma}e^{-a(x+y)}\sin(b(x+y))dydx=
$$
$$
=B_{f(\cdot)}(\alpha,\beta)\int_{0}^{\infty}g(r)r^{\alpha+\beta+\gamma-1}e^{-ar}\sin(br)dr=B_{f(\cdot)}(\alpha,\beta)\Gamma_{ \sin(b\cdot)e^{\cdot(1-a)}g(\cdot) }(\alpha+\beta+\gamma).
$$

{\it Proof} is similar to Theorem 7. $\Box$

One of the more general possibilities for the practical application of our generalized two-dimensional gamma-function is provided by the following theorem.

{\it Theorem 9. }{\it Let us assume that following conditions are met:}

 {\it 1) The function} $Q(x,y)=f\left(\frac{y}{x+y}\right)g\left( x+y \right)$ {\it satisfy the conditions of at least one of Theorems 1-4;}

 {\it 2) The function} $f(x)$ {\it is expanded into a Taylor series}
\begin{gather}\label{eqttt27}
 f(x)=\sum_{n=L}^{\infty}a_{n}x^{n}, \quad x \in [0,1], \quad L \in \mathbb{N};
\end{gather}

{\it 3) Monotone sequence of natural numbers } $\{a_{n}, n\in \mathbb{N}  \}$, {\it where} $\forall n\in \mathbb{N}: a_{n}>0$, {\it satisfies the following condition}
$$
\exists N \in \mathbb{N} \hspace{1mm} \forall  n\geq N: n a_{n}<A, \quad A<+\infty.
$$

{\it Then the following formulas will be valid}
\begin{gather}\label{eq45}
B_{f(\cdot)}(\alpha,\beta)=B(\alpha,\beta)\left(\frac{\Gamma(\alpha+\beta)}{\Gamma(\alpha)}\sum_{n=L}^{\infty}\frac{\Gamma(\alpha+n)a_{n}}{\Gamma(\alpha+\beta+n)} \right), \quad \alpha>1,  \quad \beta>1;
\end{gather}
$$
\Gamma_{\mathbb{R}^{2}_{+}; f\left(\frac{\ast_{_{\alpha}}}{\cdot_{_{\beta}}+\ast_{_{\alpha}}}\right)g\left( \cdot_{_{\beta}}+\ast_{_{\alpha}} \right)}(\alpha,\beta;\gamma)=\Gamma(\beta)\Gamma_{g(\cdot)}(\alpha+\beta+\gamma)\sum_{n=L}^{\infty}\frac{\Gamma(\alpha+n)a_{n}}{\Gamma(\alpha+\beta+n)},\quad \alpha>1, \hspace{1mm} \beta>1, \hspace{1mm} \gamma\geq 0;
$$
$$
\Gamma_{\mathbb{R}^{2}_{+}; f\left(\frac{\ast_{_{\alpha}}}{\cdot_{_{\beta}}+\ast_{_{\alpha}}}\right)g\left( \cdot_{_{\beta}}+\ast_{_{\alpha}} \right)}(\alpha,\beta;\gamma)=\Gamma_{g(\cdot)}(\alpha+\beta+\gamma)\sum_{n=L}^{\infty}\frac{a_{n}}{\Gamma_{g(\cdot)}(\alpha+\beta+\gamma+n)}\Gamma_{\mathbb{R}^{2}_{+}; g\left( \cdot_{_{\beta}}+\ast_{_{\alpha}} \right)}(\alpha+n,\beta;\gamma),
$$
$$
\quad \alpha>1, \hspace{1mm} \beta>1, \hspace{1mm} \gamma\geq 0.
$$
{\it Proof.} Multiplying both sides of equality \eqref{eqttt27} by $x^{\alpha-1}(1-x)^{\beta-1}$ and integrating from 0 to 1, we obtain
$$
B_{f(\cdot)}(\alpha,\beta)= \sum_{n=L}^{\infty}a_{n}B(\alpha+n,\beta)=\Gamma(\beta)\sum_{n=L}^{\infty}\frac{\Gamma(\alpha+n)a_{n}}{\Gamma(\alpha+\beta+n)}=
$$
$$
=B(\alpha,\beta)\left(\frac{\Gamma(\alpha+\beta)}{\Gamma(\alpha)}\sum_{n=L}^{\infty}\frac{\Gamma(\alpha+n)a_{n}}{\Gamma(\alpha+\beta+n)} \right), \quad \alpha>1,  \quad \beta>1.
$$

We obtain the last two formulas using the definition of the two-dimensional generalized gamma function and generalized formula \eqref{eq12} for the relationship between the gamma function and the beta function. $\Box$

{\it Remark 6.} Obviously, all theorems, lemmas and remarks of this section will be valid if we consider the case $Q_{2}(x,y):=f_{1}\left(\frac{y}{x+y}\right)f_{2}\left(\frac{x}{x+y}\right)g\left( x+y \right)$ of the function $Q(x,y)$. Namely, if we transform  $Q_{2}(x,y)$ in the form $Q_{2}(x,y)=f_{1}\left(\frac{y}{x+y}\right)f_{2}\left(1-\frac{y}{x+y}\right)g\left( x+y \right)$, the we can easily extend all the result to this case.

As an additional application, we can consider the hypergeometric function of the form \cite {3}:
\begin{gather}\label{eq251}
F(a,b;c;z)=\frac{\Gamma(c)}{\Gamma(b)\Gamma(c-b)}\int_{0}^{1}\frac{t^{b-1}(1-t)^{c-b-1}}{(1-tz)^{a}}dt=\frac{\Gamma(c)}{\Gamma(b)\Gamma(c-b)}B_{(1-\cdot z)^{-a}}(b,c-b),
\end{gather}
where $Re$ $(c)$ $>$ $Re$ $(b)>0$, $z$ is not a real number that is greater than or equal to $1$.

We can similarly represent this function as a special case of a two-dimensional generalized gamma function, in the form
$$
F(a,b;c;z)=\frac{1}{\Gamma(b)\Gamma(c-b)}\int_{0}^{\infty}\int_{0}^{\infty}(x+y(1-z))^{-a}(x+y)^{a}y^{b-1}x^{c-b-1}e^{-x-y}dydx=
$$
$$
=\frac{1}{\Gamma(b)\Gamma(c-b)}\Gamma_{\mathbb{R}^{2}_{+}; \left(\cdot_{_{c-b}}+\ast_{_{b}}(1-z) \right)^{-a}\left(\cdot_{_{c-b}}+\ast_{_{b}}\right)^{a}}(b,c-b;0)
$$
for the same parameters, using the equality \eqref{eq251}.

Also, if we consider a similar generalization for the hypergeometric function in the form
\begin{gather}\label{eq252}
F_{f(\cdot)}(a,b;c;z):=\frac{\Gamma(c)}{\Gamma(b)\Gamma(c-b)}\int_{0}^{1}\frac{t^{b-1}(1-t)^{c-b-1}}{(1-tz)^{a}}f(t)dt,
\end{gather}
where $Re(c)>Re(b)>0$, $z$ is not a real number that is greater than or equal to $1$, the function
$$Q_{3}(x,y):=f\left(\frac{y}{x+y}\right)g\left( x+y \right)(x+y(1-z))^{-a}(x+y)^{a}$$
it satisfy the conditions of function $Q(x,y)$ of at least one of Theorems 1-4, taking into account the remark 5, we obtain
$$
F_{f(\cdot)}(a,b;c;z)=\frac{\Gamma(c)}{\Gamma(b)\Gamma(c-b)\Gamma_{g(\cdot)}(c+\gamma)}\cdot
$$
$$
\cdot \int_{0}^{\infty}\int_{0}^{\infty}f\left(\frac{y}{x+y}\right)g\left( x+y \right)(x+y(1-z))^{-a}y^{b-1}x^{c-b-1}(x+y)^{a+\gamma}e^{-x-y}dydx=
$$
\begin{gather}\label{eq253}
=\frac{\Gamma(c)}{\Gamma(b)\Gamma(c-b)\Gamma_{g(\cdot)}(c+\gamma)}\Gamma_{\mathbb{R}^{2}_{+}; \left(\cdot_{_{c-b}}+\ast_{_{b}}(1-z) \right)^{-a}\left(\cdot_{_{c-b}}+\ast_{_{b}}\right)^{a}f\left(\frac{\ast_{_{b}}}{\cdot_{_{c-b}}+\ast_{_{b}}}\right)g\left( \cdot_{_{c-b}}+\ast_{_{b}} \right)}(b,c-b;\gamma),
\end{gather}
where $Re(c)>Re(b)>0$, $Re(\gamma)\geq 0$, $z$ is not a real number that is greater than or equal to $1$, the function $g$ identically not equal to zero.

If we assume $g\equiv 1$, and the function $Q_{4}(x,y):=f\left(\frac{y}{x+y}\right)(x+y(1-z))^{-a}(x+y)^{a}$  it satisfy the conditions of function $Q(x,y)$ of at least one of Theorems 1-4, we obtain
$$
F_{f(\cdot)}(a,b;c;z)=\frac{\Gamma(c)}{\Gamma(b)\Gamma(c-b)\Gamma(c+\gamma)}\cdot
$$
$$
\cdot \int_{0}^{\infty}\int_{0}^{\infty}f\left(\frac{y}{x+y}\right)(x+y(1-z))^{-a}y^{b-1}x^{c-b-1}(x+y)^{a+\gamma}e^{-x-y}dydx=
$$
\begin{gather}\label{eq254}
=\frac{\Gamma(c)}{\Gamma(b)\Gamma(c-b)\Gamma(c+\gamma)}\Gamma_{\mathbb{R}^{2}_{+}; \left(\cdot_{_{c-b}}+\ast_{_{b}}(1-z) \right)^{-a}\left(\cdot_{_{c-b}}+\ast_{_{b}}\right)^{a}f\left(\frac{\ast_{_{b}}}{\cdot_{_{c-b}}+\ast_{_{b}}}\right) }(b,c-b;\gamma),
\end{gather}
where $Re(c)>Re(b)>0$, $Re(\gamma)\geq 0$, $z$ is not a real number greater than or equal to $1$.

Obviously, formulas \eqref{eq253} and \eqref{eq254} for the generalized hypergeometric function \eqref{eq252} give us many options for choosing function $g$ and parameter $\gamma$.

\medskip

\textbf{Conclusion.} In the article we showed that a special case \eqref{eq11} of the two-dimensional generalized gamma function \eqref{eq1dd3d4} is the product of the one-dimensional generalized beta function \eqref{eq2} and the one-dimensional generalized gamma function \eqref{eq1}. Thus, formula  \eqref{eq12} for this special case \eqref{eq11} is a generalization of formula for the relationship between the classical gamma function and beta function. We also obtained quite a few properties of this generalization of the two-dimensional gamma function, including a formula for its
differentiation of any order. Next, we showed some practical applications of this generalization, including its use for transforming a one-dimensional generalized hypergeometric function into a two-dimensional generalized gamma function.


\medskip

\begin{thebibliography}{200}
\addcontentsline{toc}{section}{References}
\itemsep=0pt

\bibitem{1}
Andrews L.G. Special Function for Engineers and Applied Mathematics, New York: Macmillan, (1985)
\bibitem{2}
Andrews G., Askey R., Roy R. Special Function, New York: Cambridge Univ. Press, (1999)
\bibitem{3}
Bateman H., Erdelyi A. Higher Transcendental Functions, New York: McGraw-Hill, (1953) Vol. 1.
\bibitem{4}
Titchmarsh E.C. The Theory of Riemann Zeta-Function, London: Oxford Univ. Press, (1951)

\end{thebibliography}
\end{document}